%

\input amstex
\documentstyle{amsppt}
\magnification=\magstep1 
\pageheight{18cm}

 \baselineskip 24pt
\define\k{\kappa}
\define\a{\alpha}
\redefine\l{\lambda}
\define\tr{\vartriangleleft}
\define\ntr{\ntriangleleft}
\define\ra{\rightarrow}
\topmatter
\title
Any behaviour of the Mitchell Ordering of Normal Measures Is Possible\\
{\it PRELIMINARY VERSION}
\endtitle
\author  Ji\v r\' i Witzany   \endauthor
\thanks I want to express my gratitude to T.Jech
and also to J.Zapletal for many valuable discussions and 
remarks on the subject.
   \endthanks
\affil     The Pennsylvania State University  and
Charles University (Prague)                   \endaffil
\address Department of Mathematics, The Pennsylvania State University,
University Park, PA 16802                       \endaddress
\subjclass  03E35, 03E55                    \endsubjclass
\keywords        Stationary sets, reflection, measurable cardinals,
repeat points               \endkeywords

\abstract Let $U_0,U_1$ be two normal measures on $\k .$
We say that $U_0$ is in the Mitchell oredering less then $U_1,$
$U_0\tr U_1,$ if $U_0 \in Ult(V,U_1) .$
The ordering is well-known to be transitive and well-founded.
It has been an open problem to find a model where $\tr$ embeds
the four-element poset $|\; | .$ We find a generic extension
where all well-founded posets are embeddable. 
Hence there is no structural restriction on the Mitchell
ordering.
Moreover we show that it is possible to have two $\tr$-incomparable
measures that extend in a generic extension into two
$\tr$-comparable measures.

\endabstract 
\date December 1, 1993 \enddate
\email witzany\@math.psu.edu                    \endemail 
\endtopmatter

\def\lheadline{\folio\hfil {\eightpoint JI\v R\' I WITZANY}\hfil}
\def\rheadline{\hfil{\eightpoint ON THE MITCHELL ORDERING} \hfil\folio}
\headline={\ifodd\pageno\rheadline \else\lheadline\fi}

\document

We say that a well-founded poset $P$ embeds into the Mitchell ordering
of normal measures on $\k$ if there are measures $\{U_p;p\in P\}$
so that $U_p \tr U_q$ iff $p<_P q .$ 
In the well-known Mitchell's model $L[\overrightarrow U ]$
the ordering is linear (cf.[Mi83]).
S.Baldwin has constructed a model where $\tr$ is
a prewellordering (cf.[Ba85]).
 Recently
J. Cummings has described the Mitchell ordering in a particular
generic extension where it
embeds any well-founded poset that does not embed
the four element poset $|\;|$ (cf.[Cu93]). It is actually an open
problem  of [Ba84] to find  a model where the
Mitchell ordering embeds the four-element poset $|\; | .$
We show that in a generic extension 
any well-founded poset is embeddable. More
specifically:

\proclaim{Theorem}
Assume that $V$ satisfies GCH, $\k$ is measurable. Then there is
a generic extension $V[G]$ preserving cardinals, cofinalities and GCH
such that any well-foun-ded 
$\k^+$-like poset $P\in V[G]$
(i.e. $|P|\leq\k^+$ and $\forall p\in P: |P\restriction p|\leq\k$ )
such that $o(P)<o(\k)$
is embeddable into $\tr$ on normal measures over $\k$
in $V[G] .$ Moreover if $\k$ is $\Cal P_2 \k$-strong in $V$, then
any well-founded
poset of cardinality $\leq \k^+$ is embeddable.
\endproclaim

\demo{Proof}
The forcing $P_\k$ is an Easton support iteration of 
$Add(1,\l^+)^{V(P_\l)}$ ($\l<\k$ inaccessible).
That preserves cardinals, cofinalities and GCH. Let $G$ be
$P_\k$-generic/$V .$

\proclaim{Lemma 1}
Let $U_0$ be a measure on $\k,$ $j_0 : V \ra M_0 .$
Then $j_0$ can be lifted to $j_0^*:V[G] \ra M_0[G\ast \tilde H_0]$
uniquelly determined by $\tilde H_0 \in V[G] .$
\endproclaim

\demo{Proof} $j_0(P_\k)=P_\k \ast R,$ cardinality of the set
of $R/G$-antichains $\in M_0[G]$ is $\k^+$ computed in $V[G],$
and the forcing is $\k$-closed.
\qed
\enddemo

Consider only one-step extensions of this type.
Factor $\tilde H_0 = H_0^\k \ast H_0 ,$ where $H_0^\k$ is
$Add(1,\k^+)^{V[G]}$-generic/$M_0[G] .$
We want to find some sufficient and necessary conditions on $\tr$ 
on those extensions.

\proclaim{Lemma 2}
Let $U_0^*, U_1^*$ be extensions of $U_0, U_1$ given by
$H_0^\k \ast H_0$ and $H_1^\k \ast H_1 .$ If
$U_0^* \tr U_1^*$ then $U_1 \neq U_0, U_1\ntr U_0$ and 
$H_0^\k \ast H_0 \in M_1[G\ast H_1^\k] .$
On the other hand if $U_0\tr U_1$ and 
$H_0^\k \ast H_0 \in M_1[G\ast H_1^\k]$
then $U_0^* \tr U_1^* .$
\endproclaim

\demo{Proof}
Assume $U_0 \tr U_1$ and 
$H_0^\k \ast H_0 \in M_1[G\ast H_1^\k].$
Then extend $j_0:M_1 \ra M_0^1$ to 
$\tilde j_0: M_1[G] \ra M_0^1[G \ast H_0^\k \ast H_0 ]$
in $M_1[G\ast H_1^\k] .$ That defines $U_0^*$ in $M_1[G\ast H_1^\k]$
since subsets of $\k$ are same in $V[G]$ and $M_1[G] .$

Assume that $U_0^* \tr U_1^*$ then $j_0(\k)< j_1(\k).$
$U_0^* \in M_1[G\ast H_1^\k]$ since the rest of the forcing is sufficiently
closed. Consequently $j_0^*(G)=G\ast H_0^\k\ast H_0 \in M_1[G\ast H_1^\k].$
\qed
\enddemo

Now given P we are going to find $U_p^*$ ($p\in P$) such that
$U_p^* \tr U_q^*$ iff $p <_P q .$ 
Fix $U_0 \tr \cdots \tr U_\a \tr \cdots $ ($\a < o(P)$) in $V.$
Let $j_\a:V \ra M_\a$ be the corresponding embeddings

\proclaim{Claim 3}
There is $\tilde H \in V[G]$ simultaneously 
$Add(1,\k^+)^{V[G]}$-generic/$M_\a[G]$ for all $\a<o(P) .$
\endproclaim

\demo{Proof}
The cardinality of $Add(1,\k^+)^{V[G]}$-antichains $\in M_\a[G]$
computed in $V[G]$ is $\k^+$ for a fixed $\a,$
hence still $\k^+$ for all $\a<o(P)<\k^{++}$ together. $Add(1,\k^+)$
is $\k$-closed.
\qed
\enddemo

Assume that $P=(\theta, <_P),$ so that the ordering of ordinals extends $<_P .$
Factor $\tilde H = \prod_{p \in P} g_p$ using
a canonical isomorphism $Add(1,\k^+) \cong Add(\theta, \k^+) .$
For $q\in P$ define $H_q^\k =\prod_{p\leq _P q} g_p$
an $Add(1,\k^+)^{V[G]}$-generic/all $M_\a[G]$ again using
a canonical isomorphism $Add(1,\k^+)\cong Add(o.t.\{p;p\leq_P q\},\k^+) .$

\proclaim{Claim 4}
$g_p \in M_\a[G\ast H_q^\k]$ iff $p\leq_P q .$
\endproclaim

\demo{Proof}
If $p\nleq q $ then $g_p$ is generic/$M_\a[G\ast H_q^\k] ,$
hence is not in this model. The other implication is obvious.
\qed
\enddemo

$U_p^*$ will be an extension of $U_\a$ where $\a=o_P(p) .$
We have defined $H_p^\k$ and need to find an appropriate
$H_p .$
Consider $j_\a:M_{\a+1} \ra M_\a^{\a+1}$ and find
$H_p\in M_{\a+1} [G\ast H_p^\k]$
$(j_\a P_\k )^{>\k+}$-generic/$M_\a^{\a+1}[G\ast H_p^\k]$
as in lemma 1. Observe that $H_p$ is generic/$M_\a[G\ast H_p^\k]$
as well.
That defines $U_p^*$ in $M_{\a+1}[G\ast H_p^\k ] .$

\proclaim{Claim 5}
$U_p^* \tr U_q^*$ iff $p<_P q .$
\endproclaim

\demo{Proof}
Let $U_p^* \tr U_q^*$ then by lemma 1 $\a = o_P(p) < o_P(q)=\beta$
and $H_p^\k \in M_\beta [G\ast H_q^\k ]$ (e.g. using fact 6 proven bellow).
Moreover $p \leq_P q,$ hence $p<_P q ,$ since $g_p$ can be decoded
from $H_p^\k .$

Let $p<_P q,$ $\a=o_P(p) < o_P(q) =\beta .$
All we need to prove is that
$H_p^\k \ast H_p \in M_\beta [G\ast H_q^\k] .$
But $H_p^\k$ can be decoded from $H_q^\k$
(all $\k$-sequences of ordinals that are in $V[G]$ are in $M_\beta[G]$)
hence $H_p^\k \in M_\beta [G\ast H_q^\k] .$
By the construction 
$H_p\in M_{\a+1}[G\ast H _p^\k]$ and so is  
in $M_\beta[ G\ast H_p^\k ] \subset M_\beta[G\ast H_q^\k]$
using the following

\proclaim{Fact 6}
Let $U_0 \tr U_1$ then $V_{j_0(\k)+1} \cap M_0 \subset M_1 .$
\endproclaim

\demo{Proof}
Since $U_0 \in M_1$ and $V_{\k+1} \subseteq M_1$
applying $j_0$ we get
$V_{j_0(\k)+1} \cap M_0 \subseteq M_0^1 \subset M_1 .$\newline
\qed Claim 5
\enddemo
\enddemo

Finaly let $\k$ be $\Cal P_2(\k)$-strong and $P\in V[G]$
an arbitrary well-founded poset of cardinality $\kappa^+ .$
Define the generics $H_p^\k, H_p$ for $p\in P$ exactly
as above. The only problem is in the proof
of claim 5 that in general we cannot
decode $H_p^\k$ from $H_q^\k$ if $p<_P q .$
To do that we need the sequence of ordinals $\subset o.t.\{p';p'\leq_P q\}$
corresponding to the set $\{p';p'\leq p\}$ that may have cardinality
$\k^+ .$
However in that case we can assume that $P\in M_\a[G]$
for all $\a<o(P) .$ Just use a $P_\k$-name $\dot P$ for $P$ of cardinality
$\k^+$
and the Laver's diamond to get $U_0\tr \cdots \tr U_\a\tr \cdots$
such that $\dot P \in M_\a$ for all $\a<o(P) .$
\qed {Theorem}
\enddemo

\subhead Remark 1 \endsubhead
It is not true in general that $U_0^* \tr U_1^*$ implies
$U_0\tr U_1$ as lemma 2 might suggest.

Instead of the Easton iteration we could as well use
the Easton product $\tilde P_\k$ of $Add(1,\l^+)$
($\l<\k$ inaccessible).
Start with two measures $U_0\tr U_1$ in $V ,$
with the coresponding canonical embeddings $j_0:V\ra M_0, \; j_1:V\ra M_1 .$
Let $G\times \tilde G$ be $\tilde P_\k \times \tilde P_\k$-generic/$V .$
Then find  $H_0^\k\times H_1^\k \times H_2^\k \in V$ a filter 
$Add(3,\k^+)$-generic/$M_1 $ and
$H_1\times \tilde H_1 \in V$ a filter 
$(j_1\tilde P_\k)^{>\k^+}\times (j_1\tilde P_\k)^{>\k^+}$-generic/$M_1,$
and 
$H_0\times \tilde H_0 \in M_1$ a filter 
$(j_0\tilde P_\k)^{>\k^+}\times (j_0\tilde P_\k)^{>\k^+}$-generic/$M_0.$
Using Easton's lemma $G,\tilde G,H_0^\k,H_1^\k,H_2^\k,H_1,\tilde H_1$
are mutually generic/$M_1$ and 
$G,\tilde G,H_0^\k,H_1^\k,H_2^\k,H_0,\tilde H_0$ mutually generic/$M_0 .$
Firstly extend $U_0, U_1$ to $U_0^*, U_1^*$ in $V[G]$ so that
$j_1^*(G)=G\times H_1^\k \times H_1$ and 
$j_0^*(G)=G\times (H_0^\k \otimes H_1^\k)\times H_0$
where $H_0^\k \otimes H_1^\k$ denotes a coding of $H_0^\k,H_1^\k$
into  an $Add(1,\k^+)$-generic.
Obviously $U_0^* \ntr U_1^*$ since $H_0^\k \notin M_1[j_1^*(G)] .$
Then extend $U_0^*, U_1^*$ to $U_0^{**}, U_1^{**}$ in
$V[G\times \tilde G]$ so that
$j_1^{**}(\tilde G)=\tilde G \times (H_0^\k \otimes H_2^\k) \times \tilde H_1$
and
$j_0^{**}(\tilde G)=\tilde G \times H_2^\k \times \tilde H_0 .$
Then $U_0^{**} \tr U_1^{**}$ since
$U_0^*,j_0^{**}(\tilde G) \in M_1[j_1^{**}(G\times\tilde G)] .$

I do not know whethet the same can go through for the Easton iteration.

\subhead Remark 2 \endsubhead
We could use as well the Kunen--Paris forcing: the Easton product of
$Add(1,\l)$ where $\l<\k$ is a successor cardinal.

\subhead Remark 3 \endsubhead
If $V=K[\overrightarrow U_{max}],$
the core model for a coherent sequence of measures,
 where $\k$ is maximal measurable,
then we can use the method of [Cu93] to classify all measures in $V[G]:$

A finite normal iteration $j:V\ra N$ of length $n+1$ is an iteration
of ultraproducts by measures on $\k=\k_0<\k_1<\cdots\k_n .$
Any finite normal iteration $j:V\ra N$ that starts with a measure $U$
gives $\k^{++}$ extensions $U^*$ in $V[G]$ of $U$ such that 
$j_{U^*}\restriction V = j ,$ $j_{U^*}(G) = G\ast H^\k \ast H .$
And all measures in $V[G]$ are produced in this way. We can give sufficient and
necessary conditions for $\tr$ in $V[G] :$
Let $U_0^*, U_1^*$ extending $U_0, U_1$ are given by finite normal iterations
$j_0:V\ra N_0,$ $j_1:V\ra N_1$ and $H_0^\k\ast H_0,$ $H_1^\k \ast H_1 .$
Then $U_0^* \tr U_1^*$ iff $j_0\restriction Ult(V, U_1)$
is an internal iteration in this model and 
$H_0^\k \ast H_0 \in Ult(V,U_1)[G\ast H_1^\k] .$

However we can hardly describe the ordering $\tr$ in $V[G]$
in a simple manner.

That is illustrated by the following:
Let $U_0$ be the minimal measure in $V,$ let $U_0^*$ be its one-step
extension using $H_0^\k \ast H_0 \in V[G] .$
We have seen that there may be measures above $U_0^*$ even if
$H_0^\k\ast H_0 \notin M_\a[G]$ for all $\a < o(\k) .$
However it is also possible that there are no measures
above $U_0^* .$ It follows from the following joint
lemma with J. Zapletal.

\proclaim {Lemma} There is $H_0^\k \in V[G]$ $Add(1,\k^+)$-generic/$M_0[G]$
such that there is no $H_\a^\k \in V[G]$ $Add(1,\k^+)$-generic/$M_\a[G],$
$\a<o(\k)$ satisfying $H_0^\k \in M_\a[G][H_\a^\k] .$
\endproclaim

\demo{Proof}
Let $R\subseteq \k^+ \times \k^+$ be a well ordering of order type
$\gamma$ where $\gamma > \k^{++M_\a}$ for all $\a<o(\k)<\k^{++} .$
Notice that if $H_\a^\k$ is any $Add(1,\k^+)$-generic/$M_\a[G]$
then still $R \notin M_\a[G][H_\a^\k] ,$
otherwise  $\gamma$ would be less
then $\k^{++M_\a} .$
So it would be enough to code $R$ into $H_0^\k .$
Let $\langle a_\a ; \a<\k^+ \rangle \in M_1[G]$ canonically enumerate $\k^+\times\k^+$
and $\langle D_\a; \a<\k^+ \rangle \in M_1[G]$
enumerate all $M_0[G]$-dense subsets of $Add(1,\k^+),$
each set $D_\a$ enumerated by ordinals $<\k^+ .$
Construct a descending sequence of conditions
$\langle p_\a ; \a<\k^+ \rangle \in V[G]$ as follows:
assume $\langle p_\delta ; \delta < \a \rangle$ has been constructed,
then find the first $q\in D_\a$ extending $\cup\{p_\delta;\delta<\a\},$
let $\eta = \sup\{\xi +1; \xi \in dom (q)\}$
 and put $p\restriction \eta =q$ and
$$ p(\eta)=\cases 1 &\text{ iff $a_\a \in R$}\\
                          0 &\text{otherwise.} \endcases
$$
That gives an $Add(1,\k^+)$-generic filter $H_0^\k$ over $M_0[G]$ such that
$H_0^\k \in M_\a[G][H_\a^\k]$ implies $R  \in M_\a[G][H_\a^\k].$
\qed
\enddemo
 
\subhead Remark 4 \endsubhead
We can still ask what well-founded $\k^{++}$-like posets are
embeddable. For example can we embed the poset consisting of a chain
of length $\k^{++}$ and one incomparable element? Using the ideas above
we can even prove that in a sense any poset of cardinality $\k^{++}$
is embeddable.

We say that {\it a set of measures $\Cal S$ covers $\Cal P (\k^+)$}
if
$$\forall A\subseteq \k^{++} \exists U\in \Cal S :\;
   A\in Ult(V,U) .$$
If $\k$ is $\Cal P_2 \k$-strong then the measures on $\k$
cover $\Cal P(\k^+) .$ However we show that this is a much weaker
property than $\Cal P_2(\k)$-strongness, actually equiconsistent
to $o(\k) = \k^{++} .$

\proclaim{Lemma}(+Zapletal) Let $V=L[\overrightarrow U],$
$o(\k)=\k^{++} ,$ $\k$ maximal measurable. Then the measures
on $\k$ cover $\Cal P(\k^+) .$
\endproclaim

\demo{Proof}
Let $A\subseteq \k^+,$ then there is $\theta < \k^{++}$ such that
$A\in L_\theta[\overrightarrow U].$
Obviously $\overrightarrow U \cap L_\theta[\overrightarrow U]
\subseteq \overrightarrow U \restriction (\k, \theta),$
hence
$$L_\theta[\overrightarrow U] \subseteq 
L_\theta[\overrightarrow U\restriction (\k,\theta)] \subseteq
L[j_\theta(\overrightarrow U)]=
Ult(L[\overrightarrow U], U^\k_\theta) .$$
\qed
\enddemo

On the other hand we can easily destroy the covering property:
Let $\Cal S$ covers $\Cal P(\k^+) .$ Let $H$ be $Add(1,\k^+)$-generic/$V.$
Then any measure $U\in \Cal S$ remains a measure in $V[H]$
since no subsets of $\k$ have been added. On the other hand
$Ult(V[H],U)=M_U[H^*]$ is a generic extension of $M_U$
by the forcing $j_U(Add(1,\k^+))$ which is closed in the model $M_U$
-- no new subsets of $\k^+$ are added. Consequently
$H\notin Ult(V[H],U),$ and $\Cal S$ doesnot cover $\Cal P(\k^+)$
in $V[H].$

Using this observation we can prove:

\proclaim{Proposition} Assume that measures on $\k$ in $V$ cover
$\Cal P(\k^+) .$ Then there is a generic extension $V[G]$ of $V$
such that any well-founded $\k^{++}$-like poset in $V$ is embeddable
into $\tr$ in $V[G] .$
\endproclaim

\demo{Proof}
Let  $G=H\times \tilde G$ be 
$Add(\k^{++},\k^+)\times P_\k$-generic/$V$ where $P_\k$ is the Kunen-Paris
forcing. Observe that if $P$ is $\k^{++}$-like then $P$ can be enumerated
as $P=\langle p_\a;\a<\k^{++} \rangle$ so that the well-ordering
extends $<_P .$
Consequently we can assume that $P=\langle \k^{++}, <_P \rangle .$
Using the covering property in $V$ find 
$U_0\tr \cdots \tr U_\a \tr \cdots$ ($\a<\k^{++}$) such that
$\{ \gamma;\gamma <_P \a \} \in M_\a=Ult(V,U_\a) .$
Factor $H$ as $\prod _{\a<\k^{++}} g_\a$ where $g_\a$
is an $Add(1,\k^+)$-generic.

Let $H_\a^\k =\prod _{\gamma \leq _P \a} g_\gamma$ be 
an $Add(1,\k^+)$-generic using an appropriate canonical isomorphism.
Let $H_\a \in M_{\a+1}$ be $(j_\a P_\k)^{>\k}$-generic/$M_\a .$
Then $\tilde G \times H_\a^\k \times H_\a $ is
$j_\a P_\k$-generic/$Ult(V[H], U_\a)$ defining an extension $U_\a^*$
of $U_\a $ in $V[G] .$
Since $H_\a \times H_\a^\k \in \tilde M_\beta [\tilde G \times H_\beta ^\k]$
iff $\a<_P \beta$ we are done: $U_\a^* \tr U_\beta^*$ iff $\a<_P \beta .$
\qed
\enddemo

\enddocument
\end